\documentclass[11pt]{article}

\usepackage[applemac]{inputenc} 
\usepackage{amsthm,amssymb,amsbsy,amsmath,amsfonts,amssymb,amscd,mathrsfs}

\long\def\symbolfootnote[#1]#2{\begingroup%
\def\thefootnote{\fnsymbol{footnote}}\footnote[#1]{#2}\endgroup} 

\usepackage{graphicx}
\DeclareGraphicsExtensions{.jpg}
\usepackage{color}
\usepackage{fullpage}
\usepackage[hypertex]{hyperref}
\theoremstyle{plain}
\newtheorem{theorem}{Theorem}
\newtheorem{corollary}[theorem]{Corollary}
\newtheorem{lemma}[theorem]{Lemma}
\newtheorem{proposition}[theorem]{Proposition}
\theoremstyle{definition} 
\newtheorem{definition}[theorem]{Definition}
\newtheorem{remark}[theorem]{Remark}
\newcommand{\bt}{\begin{theorem}}
\newcommand{\et}{\end{theorem}}
\newcommand{\bl}{\begin{lemma}}
\newcommand{\el}{\end{lemma}}
\newcommand{\bp}{\begin{proposition}}
\newcommand{\ep}{\end{proposition}}
\newcommand{\bc}{\begin{corollary}}
\newcommand{\ec}{\end{corollary}}
\newcommand{\bdeff}{\begin{definition}}
\newcommand{\edeff}{\end{definition}}
\newcommand{\brem}{\begin{remark}}
\newcommand{\erem}{\end{remark}}

\newcommand{\bi}{\begin{itemize}}
\newcommand{\iii}{\item}
\newcommand{\ei}{\end{itemize}}
\newcommand{\bd}{\begin{description}}
\newcommand{\ed}{\end{description}}

\newcommand{\bqn}{\begin{eqnarray}}
\newcommand{\eqn}{\end{eqnarray}}
\newcommand{\eqnn}{\nonumber\end{eqnarray}}

\newcommand{\nn}{\nonumber}
\newcommand{\ba}[1]{\begin{array}{#1}}
\newcommand{\ea}{\end{array}}
\newcommand{\wh}[1]{\widehat{#1}}
\newcommand{\wt}[1]{\widetilde{#1}}
\newcommand{\lam}{\lambda}

\newcommand{\g}{\gamma}
\newcommand{\al}{\alpha}
\newcommand{\eps}{\varepsilon}

\newcommand{\ph}{\varphi}

\newcommand{\R}{\mathbb{R}}
\newcommand{\N}{\mathbb{N}}

\newcommand{\mb}[1]{\mathbb{ #1 }}
\newcommand{\mc}[1]{\mathcal{ #1 }}

\newcommand{\all}{\forall\,}

\newcommand{\la}{\langle}
\newcommand{\ra}{\rangle}

\newcommand{\VecM}{\mathrm{Vec}(M)}
\newcommand{\virg}[1]{``#1''}
\newcommand{\tx}[1]{\mathrm{#1}}
\newcommand{\til}[1]{\widetilde{#1}}

\newcommand{\der}[2]{\frac{d}{d#1}\Big|_{#1=#2}}

\newcommand{\distr}{\mc{D}}
\newcommand{\dil}{\delta}
\newcommand{\metr}{\textsl{g}}
\newcommand{\sub}{\mathbf{S}}
\newcommand{\Pg}[1]{\left\{ #1 \right\}}

\newcommand{\hp}{hypothesis}

\newcommand{\lapl}{\Delta}
\newcommand{\dive}{\text{div}}
\newcommand{\grad}{\nabla}
\newcommand{\Y}{\mc{Y}}

\begin{document}
\begin{center} \noindent
{\LARGE{\sl{\bf Trace heat kernel asymptotics in 3D contact\\[0.2cm] sub-Riemannian geometry}}}
\vskip 0.6 cm
Davide Barilari\\ 
{\footnotesize SISSA, Via Bonomea 265, Trieste, Italy - {\tt davide.barilari@sissa.it}}\\
\vskip 0.3cm

\vskip 0.3cm
\end{center}
\vskip 0.1 cm

\begin{center}
\today
\end{center}
\vskip 0.2 cm

\begin{abstract}
In this paper we 
study the small time asymptotics for the heat kernel on a sub-Riemannian manifold, using a perturbative approach. We  explicitly compute, in the case of a 3D contact structure, the first two coefficients of the small time asymptotics expansion of the heat kernel on the diagonal, expressing them in terms of the two basic functional invariants $\chi$ and $\kappa$ defined on a 3D contact structure.
\end{abstract}


\section{Introduction}

In this paper we consider the problem of finding the small time asymptotic expansion on the diagonal of the heat kernel on a sub-Riemannian manifold $M$.\symbolfootnote[0]{This research has been supported by the European Research Council, ERC StG 2009 \virg{GeCoMethods}, contract number 239748, by the ANR Project GCM, program \virg{Blanche}, project number NT09-504490.}

By  a sub-Riemannian manifold we mean a triple $\sub=(M,\distr,\metr)$, where $M$ is a connected orientable smooth manifold of dimension $n$, $\distr$ is a smooth vector distribution of constant rank $k<n$, satisfying the H\"ormander condition and $\metr$ is an Euclidean structure on $\distr$. 

It is well known that a sub-Riemannian manifold has a natural structure of metric space, where the distance is the so called Carnot-Caratheodory distance
\begin{eqnarray}
d(q_0,q_1)=
\inf\{\int_0^T\sqrt{\metr_{\g(t)}(\dot\g(t),\dot\g(t))}~dt~|~ \g:[0,T]\to M \mbox{ is a Lipschitz curve},\cr
\g(0)=q_0,\g(T)=q_1, ~~\dot \g(t)\in\distr_{\g(t)}\mbox{ a.e. in $[0,T]$} \}.\nonumber
\end{eqnarray}
As a consequence of the H\"ormander condition, this distance is always finite and continuous, and induces on $M$ the original topology (Chow Theorem, see e.g. \cite{agrachevbook}).

Define inductively the vector distributions $\distr^{1}:=\distr, \distr^{i+1}:=\distr^{i}+[\distr^{i},\distr]$, for every $i=1,2,\ldots$ and denote by $\distr^{i}_{q}$ the fiber of $\distr^{i}$ at the point $q$. Under the hypothesis that the sub-Riemannian manifold is regular, i.e.  for every $  i=1,\ldots,m$, the dimension of $\distr^{i}_{q}$ does not depend on the point, the H\"ormander condition guarantees that there exists a mimimal $m\in \N$, called \emph{step} of the structure, such that $\distr_{q}^{m}=T_{q}M$, for all $q\in M$. 


Moreover, it is well known that the Hausdorff dimension of $M$, as a metric space, is given by the formula (see \cite{mitchell})

$$Q=\sum_{i=1}^m i k_i, \qquad k_i:=\tx{dim}\, \distr^i - \dim \distr^{i-1}.$$
In particular the Hausdorff dimension of a sub-Riemannian manifold is always bigger than its topological dimension.

A sub-Riemannian structure on a three-dimensional manifold $M$, is said to be \emph{contact} if its
distribution $\distr$ is the kernel of a differential
one form $\distr= \ker \omega$, such that
$d\omega\wedge\omega$ is a nonvanishing
$3$-form on $M$.

The geometry of 3D contact sub-Riemannian structures have been deeply studied and they possess two basics
differential invariants $\chi$ and $\kappa$ (see \cite{agrexp} for the
precise definition and their role in the
asymptotic expansion of the sub-Riemannian exponential map). 
These invariants $\chi$ and $\kappa$ are smooth real functions on
$M$ that reflect geometric properties of the structure, and in the case of left-invariant metrics on Lie groups they have been used in \cite{miosr3d} to provide a complete classification of such structures up to local isometries.

Our goal is to compute the first term in the asymptotic expansion of the heat kernel in the 3D contact case and to understand how the invariants $\chi$ and $\kappa$ enters in the coefficients.

The heat equation on a sub-Riemannian manifold is a natural model for the description of a non-isotropic diffusion process on a manifold. It is defined by the second order PDE  
\bqn \label{eq:pde}
 \frac{\partial}{\partial t}\psi(t,x)=\lapl_{f}\psi(t,x), \qquad \all t>0,\, x\in M,
 \eqn
where $\lapl_{f}$ is the \emph{sub-Riemannian Laplacian} (also called \emph{sub-Laplacian}) which is a hypoelliptic, but not elliptic, second order differential operator. Locally this operator can be written in the form \virg{sum of squares} plus a first order part
$$\lapl_{f}=\sum_{i=1}^{k}f_{i}^{2}+a_{i}f_{i}, \qquad \qquad k=\dim \distr <n.$$
where $f_{1},\ldots,f_{k}$ is an orthonormal basis for the sub-Riemannian structure and $a_{1},\ldots,a_{k}$ are suitable smooth coefficients (see Section \ref{s:lapl} for the precise meaning of these coefficients).

From the analytical viewpoint, these operators, and their parabolic counterpart $\partial_{t}-\lapl_{f}$, have been widely studied, starting from the celebrated work of H\"ormander \cite{hormander}.
A probabilistic approach to hypoelliptic diffusion equation can be found in \cite{baudoin,bismut,kusustroock}, where the existence of a smooth heat kernel for such equations is given.

On the other hand, a \virg{geometric} definition of the Laplacian is needed if one want to 
find some relation between the analytical properties of the heat kernel (e.g. the small time asymptotics)
and the geometric properties of the manifold, like in the Riemannian case (see \cite{rosenberg,spectrebook} for the relation between the heat kernel and the Riemannian curvature of the manifold, and \cite{neel} for a characterization of the cut locus via the heat kernel).

As it was pointed out in \cite{laplacian,montgomerybook}, to have an intrinsic definition of the sub-Laplacian $\lapl_{f}$ (i.e. that depends only on the geometric structure) it is necessary to build an instrinsic volume for the structure.

In the sub-Riemannian case there are two intrinsic volumes that are defined, namely the Popp's volume (that is the analogue of the Riemannian volume form in Riemannian geometry) and the Hausdorff volume. In \cite{miovolume} it is proved that, starting from dimension 5, in general they are not proportional. On the other hand in the 3D contact case they coincide.

The existence of an asymptotic expansion for the heat kernel was proved, besides the classical Riemannian case, when the manifold is endowed with a time dependent Riemannian metric in \cite{lanconelligarofalo}, in the sub-Riemannian free case (when $n=k+\frac{k(k-1)}{2}$) in \cite{brockettmansouri}. In  \cite{leandre,benarouscut,taylor} the general sub-Riemannian case is considered, using a probabilistic approach, obtaining different expansions on the cut locus and out of that.

The same method was also applied in \cite{benarousdiag} to obtain the asymptotic expansion on the diagonal. In particular it was proved that, for the sub-Riemannian heat kernel $p(t,x,y)$ the following expansion holds
\bqn \label{eq:uno}
p(t,x,x)\sim \frac{1}{t^{Q/2}}(a_{0}+a_{1}t+a_{2}t^{2}+\ldots+a_{k}t^{k}+O(t^{k+1})), \qquad \text{for}\quad t\to 0,
\eqn
where $Q$ denotes the Hausdorff dimension of $M$.

Besides these existence results, the geometric meaning of the coefficients in the expansion on the diagonal (and out of that) is far from being understood, even in the simplest 3D case, where the heat kernel has been computed explicitly in some cases of left-invariant structures on Lie groups only (see \cite{laplacian,bonnefontsu2,taylorparam}). In analogy to the Riemannian case, one would expect that the curvature tensor of the manifold and its derivatives appear in 
these expansions.

In this paper we compute the first terms (precisely $a_{0}$ and $a_{1}$, referring to \eqref{eq:uno}), for every 3D contact structure. Our main tool is the nilpotent approximation of the sub-Riemannian structure.
Under the regularity hypothesis, the nilpotent approximation of a sub-Riemannian structure $\sub=(M,\distr,\metr)$ at a point $q\in M$, is endowed itself with a left-invariant sub-Riemannian structure $\widehat{\sub}_{q}$ on a so-called Carnot group (i.e. homogeneous nilpotent and simply connected Lie group) whose Lie algebra is generated by the nilpotent approximation of a basis of the Lie algebra of $\sub$, as explained in Section \ref{s:nilpapp}.  

When $\sub$ is a 3D contact sub-Riemannian manifold, the nilpotent approximation $\wh{\sub}_{q}$ of $\sub$ at every point $q\in M$ is isometric to the Heisenberg group, i.e. the sub-Riemannian structure on $\R^{3}$ (where coordinates are denoted by $q=(x,y,w)$) defined by the global orthonormal frame
\begin{align} \label{eq:heisenberg}
\wh{f}_{1}&=\partial_{x}+\frac{y}{2}\partial_{w}, \qquad \wh{f}_{2}=\partial_{y}-\frac{x}{2}\partial_{w}.
\end{align} 

Using this approach, we reduce the problem of computing the small time asymptotic of the heat kernel of the original sub-Riemannian structure to the problem of studying the heat kernel of a family of approximating structures at the fixed time $t=1$. In such a way we present the original structure (and as a consequence its heat kernel) as a perturbation of the nilpotent one. With the perturbative method we then compute the coefficient $a_{0}$, which reflects the properties of the Heisenberg group case (being the nilpotent approximation of every 3D contact structure), and the coefficient $a_{1}$, where the local invariant $\kappa$ appears.

The main result is stated as follows:

\bt\label{t:main} Let $M$ be a 3D contact sub-Riemannian structure, with local invariants $\chi$ and $\kappa$. Let $p(t,x,y)$ denotes the heat kernel of the sub-Riemannian heat equation. Then the following small time asymptotic expansion hold
$$p(t,x,x) \sim \frac{1}{16t^{2}}(1+\kappa(x)t+ O(t^{2})), \qquad \text{for}\quad t\to 0.$$
\et
%
\noindent
Notice that the Hausdorff dimension of a 3D contact structure is $Q=4$.

\section{Basic definitions}
We start recalling the definition of sub-Riemannian manifold.
\bdeff
A \emph{sub-Riemannian manifold} is a triple $\sub=(M,\distr,\metr)$, 
where
\bi
\iii[$(i)$] $M$ is a connected orientable smooth manifold of dimension $n\geq 3$;
\iii[$(ii)$] $\distr$ is a smooth distribution of constant rank $k< n$ satisfying the \emph{H\"ormander condition}, i.e. a smooth map that associates to $q\in M$  a $k$-dimensional subspace $\distr_{q}$ of $T_qM$ and we have
\bqn \label{Hor}
\text{span}\{[X_1,[\ldots[X_{j-1},X_j]]](q)~|~X_i\in\overline{\distr},\, j\in \N\}=T_qM, \quad \all q\in M,
\eqn
where $\overline{\distr}$ denotes the set of \emph{horizontal smooth vector fields} on $M$, i.e. $$\overline{\distr}=\Pg{X\in\mathrm{Vec}(M)\ |\ X(q)\in\distr_{q}~\ \forall~q\in M}.$$
\iii[$(iii)$] $\metr_q$ is a Riemannian metric on $\distr_{q}$ which is smooth 
as a function of $q$. We denote  the norm of a vector $v\in \distr_{q}$ with $|v|$, i.e.  $|v|=\sqrt{\metr_{q}(v,v)}.$
\ei
\edeff

A Lipschitz continuous curve $\g:[0,T]\to M$ is said to be \emph{horizontal} (or \emph{admissible}) if 
$$\dot\g(t)\in\distr_{\g(t)}\qquad \text{ for a.e. } t\in[0,T].$$

Given an horizontal curve $\g:[0,T]\to M$, the {\it length of $\g$} is
\bqn
\label{e-lunghezza}
\ell(\g)=\int_0^T |\dot{\g}(t)|~dt.
\eqn
The {\it distance} induced by the sub-Riemannian structure on $M$ is the 
function
\bqn
\label{e-dipoi}
d(q_0,q_1)=\inf \{\ell(\g)\mid \g(0)=q_0,\g(T)=q_1, \g\ \mathrm{horizontal}\}.
\eqn
The \hp\ of connectedness of $M$ and the H\"ormander condition guarantee the finiteness and the continuity of $d(\cdot,\cdot)$ with respect to the topology of $M$ (Chow-Rashevsky theorem, see for instance \cite{agrachevbook}). The function $d(\cdot,\cdot)$ is called the \emph{Carnot-Caratheodory distance} and gives to $M$ the structure of metric space (see \cite{bellaiche}).


Locally, the pair $(\distr,{\mathbf g})$ can be given by assigning a set of $k$ smooth vector fields spanning $\distr$ and that are orthonormal for ${\mathbf g}$, i.e.  
\bqn
\label{trivializable}
\distr_{q}=\text{span}\{f_1(q),\dots,f_k(q)\}, \qquad \qquad \metr_q(f_i(q),f_j(q))=\delta_{ij}.
\eqn
In this case, the set $\Pg{f_1,\ldots,f_k}$ is called a \emph{local orthonormal frame} for the sub-Riemannian structure.

\bdeff
Let $\distr$ be a distribution. Its \emph{flag} is the sequence of distributions $\distr^{1}\subset\distr^{2}\subset\ldots$ defined through the recursive formula
$$\distr^{1}:=\distr,\qquad  \distr^{i+1}:=\distr^{i}+[\distr^{i},\distr].$$
A sub-Riemannian manifold is said to be \emph{regular} if for each $i=1,2,\ldots$ the dimension of $\distr^{i}_{q}$ does not depend on the point $q\in M$.
\edeff

A sub-Riemannian manifold is said to be \emph{nilpotent} if there exists an orthonormal frame for the structure $\{f_{1},\ldots,f_{k}\}$ and $j\in \mb{N}$ such that
$[f_{i_{1}},[f_{i_{2}},\ldots,[f_{i_{j-1}},f_{i_{j}}]]]=0$, for every commutator of length $j$.

\subsection{3D contact case}
In this section we focus on three dimensional sub-Riemannian manifolds, endowed with a contact distribution. 

\bdeff Let $M$ be a smooth manifold, with $\dim M=3$. A sub-Riemannian
structure on $M$ is said to be \emph{contact} if $\distr$ is a
contact distribution, i.e. $\distr= \ker \omega$, where $\omega\in\Lambda^{1}M$ satisfies $\omega \wedge d\omega \neq 0$.\edeff

\brem \label{rem:contact} Notice that a contact structure is forced to be bracket generating. The fact that $\omega$ is a contact form can be rewritten as $d\omega|_{\distr} \neq 0$. Moreover, it is possible to normalize the contact fom $\omega$ in such a way that $d\omega|_{\distr}$ coincides with the volume form on $\distr$ naturally defined by $\metr$ (see \cite{miosr3d, agrexp}). 
\erem

Since we deal with local properties, it is not restrictive to assume that
the sub-Riemannian structure is defined by a local orthonormal frame $f_{1},f_{2}$ that satisfies:
\begin{gather}
(M,\omega) \text{ is a 3D contact structure}, \notag  \\
\distr=\tx{span}\{f_1,f_2\}=\ker \omega,\label{eq:setting}\\
\metr(f_i,f_j)=\delta_{ij}, \quad
d\omega(f_1,f_2)=1. \notag
\end{gather}
Recall that the orthonormal frame
$f_1,f_2$ is not unique. Indeed every rotated frame (where the angle
of rotation depends smoothly on the point) defines the same
structure.

The \emph{Reeb vector field} associated to the contact structure is the unique vector field $f_0$ such that
\begin{align} \label{eq:deff0}
\omega(f_0)&=1, \notag\\
d\omega(f_0,\cdot)&=0.
\end{align}
Clearly, $f_0$ depends on the sub-Riemannian structure (and its orientation) only and not on the frame selected.

The Lie algebra of vector fields generated by $f_0,f_1,f_2$ satisfies the following commutation relations
\begin{align} \label{eq:algebracampi}
[f_1,f_0]&=c_{01}^1 f_1+c_{01}^2 f_2, \notag\\
[f_2,f_0]&=c_{02}^1 f_1+c_{02}^2 f_2, \\
[f_2,f_1]&=c_{12}^1 f_1+c_{12}^2 f_2+f_0, \notag
\end{align}
where $c_{ij}^{k}$ are functions on the manifold, also called structure constants of the Lie algebra.

Next we recall the definition of the local invariants of a contact three-dimensional structure. Here we simply give their expression when computed by means of coefficients of \eqref{eq:algebracampi}. For a geometric interpretation of these invariants and an intrinsic definition one can see  \cite{agrexp,miosr3d}.
\bdeff With reference to notation \eqref{eq:algebracampi}, we define
\bi
\iii[-]
the \emph{first invariant} \bqn \label{eq:chi}
 \chi(q)=\sqrt{-\det C}, \qquad C=
 \begin{pmatrix}
 c_{01}^{1}&(c_{01}^{2}+c_{02}^{1})/2\\
 (c_{01}^{2}+c_{02}^{1})/2&c_{02}^{2}
 \end{pmatrix},
\eqn

\iii[-] the \emph{second invariant}
\bqn \label{eq:defkappa}
\kappa(q)=f_2(c_{12}^1)-f_1(c_{12}^2)-(c_{12}^1)^2-(c_{12}^2)^2+
\dfrac{c_{01}^2-c_{02}^1}{2}.
\eqn     
\ei
\edeff
\brem \label{r:chipositivo} Notice that, from their very definition, $\chi$ and $\kappa$ are smooth functions on $M$ that do not depend on the choice of the orthonormal frame. Moreover  $\chi\geq0$ and it vanishes everywhere if and only if the flow of the Reeb vector field $f_{0}$ is a flow of sub-Riemannian isometries for $M$.
\erem

\section{The sub-Laplacian in a sub-Riemannian manifold}\label{s:lapl}
In this section we compute the intrinsic hypoelliptic Laplacian on a regular sub-Riemannian manifold $(M, \distr, \metr)$, also called \emph{sub-Laplacian}, writing its expression in a local orthonormal frame. In particular we find its explicit expression in the 3D contact sub-Riemannian case in terms of \eqref{eq:algebracampi}. 

The sub-Laplacian is the natural generalization of the Laplace-Beltrami operator $\lapl$ defined on a Riemannian manifold, that is $\lapl \phi = \dive(\grad \phi)$, where $\grad$ is the unique operator from $C^{\infty}(M)$ to $\VecM$ satisfying $$\metr( \grad \phi, X) = d\phi(X), \qquad \all\,  X\in \VecM.$$ Here $\metr$ denotes the Riemannian metric, and the divergence of a vector field $X$ is the unique function $\dive\,  X$ satisfying 
\bqn \label{eq:dive} 
L_{X}\mu=(\dive\, X )\mu,
\eqn 
where $\mu$ is the Riemannian volume form and $L_{X}$ denotes the Lie derivative. 

In the sub-Riemannian case these definitions are replaced by the notions of horizontal gradient and of divergence with respect to the Popp measure, which is well defined in the regular case (see \cite{laplacian}).

\bdeff Let $M$ be a sub-Riemannian manifold and $\phi \in C^{\infty}(M)$. The \emph{horizontal gradient} of $\phi$ is the unique horizontal vector field  $\grad \phi \in \overline{\distr}$ that satisfies
\bqn \label{eq:horgrad}
\metr( \grad \phi, X) = d\phi(X), \qquad \all\,  X\in \overline{\distr}.
\eqn
\edeff
Given a local orthonormal frame $\{f_{1},\ldots,f_{k}\}$ for the sub-Riemannian structure, it is easy to see that the horizontal gradient $\grad \phi \in \overline{\distr}$ of a function  is computed as follows
\bqn \label{eq:grad}
\grad \phi= \sum_{i=1}^{k}f_{i}(\phi)f_{i},  \qquad \phi \in C^{\infty}(M),
\eqn
where the vector field acts on functions as a differential operator. 

\medskip
{\bf Notation}. In what follows we will denote by $\lapl_{f}$ the sub-Laplacian associated to the sub-Riemannian structure defined by the local orthonormal frame $f_{1},\ldots,f_{k}$. Actually this definition does not depend on the choice of the orthonormal frame (see also Proposition \ref{p:evita}).

\medskip
In the sub-Riemannian regular case, even if there is no scalar product defined in $T_{q}M$, we can still define an intrinsic volume, called Popp volume, by means of the Lie bracket of the horizontal vector fields (see \cite{laplacian}). 
Here we recall a convenient definition of Popp measure only for the 3D contact case.
\bdeff Let $M$ be an orientable 3D contact sub-Riemannian structure and $\{f_{1}, f_{2}\}$ a local orthonormal frame. Let $f_{0}$ be the Reeb vector field and $\nu_{0},\nu_{1},\nu_{2}$ the dual basis of 1-form, i.e. $\la \nu_{i}, f_{j}\ra=\delta_{ij}$. The \emph{Popp volume} on $M$ is the form\footnote{In \cite{laplacian} the dual basis of the frame $\{f_{1},f_{2},[f_{1},f_{2}]\}$ was considered to build the Popp volume. From \eqref{eq:algebracampi} it follows that these two constructions agree each other.} $\mu := \nu_{0}\wedge \nu_{1} \wedge \nu_{2}$. 
\edeff

Using the formula
\bqn\label{eq:dive}
\dive (aX)=Xa + a \, \dive X, \qquad \all a\in C^{\infty}(M),\,  X \in \VecM,
\eqn
it is easy to find the expression of the sub-Laplacian with respect to any volume
\begin{align*}
\dive(\grad \phi)&= \sum_{i=1}^{k}\dive (f_{i}(\phi)f_{i})\\
&=\sum_{i=1}^{k}f_{i}(f_{i}(\phi))+ (\dive\, f_{i}) f_{i}(\phi).
\end{align*}
Thus
\bqn \label{eq:lapldiv}
\lapl_{f}=\sum_{i=1}^{k}f_{i}^{2}+ (\dive\, f_{i}) f_{i}.
\eqn
\brem \label{r:superem} Here we collect few properties of the sub-Laplacian that immediately follows from the definition:
\bi
\iii[(i)] The sub-Laplacian is always presented as sum of squares of the horizontal vector fields plus a first order horizontal part, whose coefficients heavily depends on the choice of the volume (see \eqref{eq:lapldiv}).
Moreover $\lapl_{f}$ is the sum of squares if and only if all the vector fields of the orthonormal frame are divergence free. 

\iii[(ii)] From  \eqref{eq:dive}  and \eqref{eq:lapldiv} it easily follows that the sub-Laplacian is a homogeneous differential operator of degree two with respect to dilations of the metric structure. More precisely, if we consider the dilated structure (denoted by $\lam f$) where all vector fields of the orthonormal frame are multiplied by  a positive constant $\lam>0$, we have
\bqn \label{eq:hom2}
\lapl_{\lam f}=\lam^{2} \lapl_{f}.
\eqn

\iii[(iii)] Define for every pair of functions $\phi,\ph \in C^{\infty}_{0}(M)$ the bilinear form
$$(\phi,\ph)_{2}= \int_{M}\phi\ph \,d\mu.$$
Given a vector field $X\in \VecM$, its formal adjoint $X^{*}$ is the differential operator that satisfies the identity
$$(X\phi,\ph)_{2}= (\phi,X^{*}\ph)_{2}, \qquad \all \phi,\ph \in C^{\infty}_{0}(M).$$
It is easily computed that $X^{*}=-X-\dive\, X$. In particular it follows that the sub-Laplacian $\lapl_{f}$ is rewritten as
$$\lapl_{f}=-\sum_{i=1}^{k} f_{i}^{*}f_{i},$$
and that satisfies the identities
\bqn \label{eq:sym}
(\lapl_{f}\phi, \ph)_{2}=(\phi,\lapl_{f} \ph)_{2}, \qquad (\phi,\lapl_{f} \phi)_{2}\leq0, \qquad \all \phi,\ph \in C^{\infty}_{0}(M).
\eqn
\ei
\erem

From this we can explicitly compute the sub-Laplacian in the 3D contact case (see also \cite{laplacian}, formula (8)).
\bp \label{p:lapl} Let $M$ be a 3D sub-Riemannian manifold and $f_{1},f_{2}$ be an orthonormal frame and $f_{0}$ be the Reeb vector field. 
Then the sub-Laplacian is expressed as follows
\bqn \label{eq:lapl3d}
\lapl_{f}=f_{1}^{2}+f_{2}^{2}+c_{12}^{2}f_{1}-c_{12}^{1}f_{2},
\eqn
where $c_{12}^{1},c_{12}^{2}$ are the structure constant appearing in \eqref{eq:algebracampi}.
\ep

\begin{proof} 
Using \eqref{eq:lapldiv}, it is enough to compute the functions $a_{i}:= \dive_{\mu} f_{i}$, where $\mu = \nu_{1} \wedge \nu_{2} \wedge \nu_{0}$ is the Popp measure.

To this purpose let us compute the quantity $L_{X}\mu$, for every vector field $X$. Recall that the action of the Lie derivative on a differential 1-form  is defined as
$$L_{X} \nu= \der{t}{0}e^{tX*}\nu, \qquad \nu \in \Lambda^{1}(M).$$
where $e^{tX}$ denotes the flow on $M$  generated by the vector field $X$.
Using the fact that $L_{X}$ is a derivation we get
\begin{align} \label{eq:sasa}
L_{X}(\nu_{0}\wedge \nu_{1} \wedge \nu_{2})&=L_{X}\nu_{0} \wedge \nu_{1} \wedge \nu_{2}+\nu_{0}\wedge L_{X}\nu_{1} \wedge \nu_{2} +\nu_{0}\wedge \nu_{1} \wedge L_{X}\nu_{2}.
\end{align}

Moreover, for every $i=0,1,2,$ we can write
$$L_{X}\nu_{i}=\sum_{j=0}^{2}a_{ij} \nu_{j},$$ 
The coefficients $a_{ij}$ can be computed evaluating $L_{X}\nu_{i}$ on the dual basis
\begin{align*}
a_{ij}&=\la L_{X}\nu_{i}, f_{j}\ra\\
&=\la\der{t}{0}e^{tX*}\nu_{i},f_{j}\ra\\
&= \la\nu_{i},\der{t}{0}e^{tX}_{*}f_{j}\ra\\
&=\la \nu_{i},[f_{j},X]\ra.
\end{align*}

Plugging these coefficients into \eqref{eq:sasa} we get
$$L_{X}\mu = ( \la \nu_{1},[f_{1},X]\ra+\la \nu_{2},[f_{2},X]\ra+ \la \nu_{0},[f_{0},X]\ra) \mu,$$
from which it follows
\bqn \label{eq:divee}
\dive\, X= \la \nu_{1},[f_{1},X]\ra+\la \nu_{2},[f_{2},X]\ra+\la \nu_{0},[f_{0},X]\ra.
\eqn

Using \eqref{eq:algebracampi} it is easy to see that $\la \nu_{0},[f_{0},X]\ra=0$ for every horizontal vector field. Thus,
applying \eqref{eq:divee} with $X=f_{i},$ with  $i=1,2$ one gets   
$$\dive\, f_{1}=\la \nu_{2},[f_{2},f_{1}]\ra=c_{12}^{2}, \qquad \dive\, f_{2}=\la \nu_{1},[f_{1},f_{2}]\ra=-c_{12}^{1}.$$
Then \eqref{eq:lapl3d} easily follows from \eqref{eq:lapldiv}.
\end{proof}

From the above construction it is clear that the sub-Laplacian depends only on the sub-Riemannian structure and not on the frame selected, i.e. it is invariant for rotations of the orthonormal frame. Here we give also a direct proof of this fact.

%

\bp \label{p:evita} The sub-Laplacian is invariant with respect to rotation of the orthonormal frame.
\ep

\begin{proof} Let us consider an orthonormal frame $f_{1},f_{2}$ and the rotated one
\begin{align}
\widetilde{f}_{1}&= \cos \theta f_{1} +\sin \theta f_{2},\\
\widetilde{f}_{2}&=-\sin \theta f_{1} + \cos \theta f_{2},
\end{align}
where $\theta=\theta(q)$ is a smooth function on $M$. If we denote by $\widetilde{c}_{ij}^k$ the structure constants computed in the rotated frame, from the formula
\bqn
[\wt{f}_{1},\wt{f}_{2}]=[f_{1},f_{2}]-f_{1}(\theta)f_{1}-f_{2}(\theta)f_{2},
\eqn
it is easy to prove that the new structure constant are computed according to the formulas
\begin{align*}
\widetilde{c}_{12}^1= \cos \theta (c_{12}^1-f_1(\theta))- \sin
\theta (c_{12}^2-f_2(\theta)), \\
\wt{c}_{12}^2= \sin \theta (c_{12}^1-f_1(\theta))+ \cos
\theta (c_{12}^2-f_2(\theta)).
\end{align*}
From these relations one can easily compute
\bqn \label{eq:lapl1}
\wt{f}_{1}^{2}+
\wt{f}_{2}^{2}= f_{1}^{2}+f_{2}^{2}+f_{1}(\theta)f_{2}-f_{2}(\theta)f_{1},
\eqn
and
\bqn \label{eq:lapl2}
-\wt{c}_{12}^{1}\wt{f}_{2}+\wt{c}_{12}^{2}\wt{f}_{1}=-(c_{12}^{1}-f_{1}(\theta))f_{2}+(c_{12}^{2}-f_{2}(\theta))f_{1}.
\eqn
Combining \eqref{eq:lapl1} and \eqref{eq:lapl2} one gets, denoting $\lapl_{\til{f}}$ the laplacian defined by the rotated frame,
\begin{align*}
\lapl_{\til{f}}&=\wt{f}_{1}^{2}+\wt{f}_{2}^{2}-\wt{c}_{12}^{1}\wt{f}_{2}+\wt{c}_{12}^{2}\wt{f}_{1}\\
&
=f_{1}^{2}+f_{2}^{2}+f_{1}(\theta)f_{2}-f_{2}(\theta)f_{1}-\wt{c}_{12}^{1}\wt{f}_{2}+\wt{c}_{12}^{2}\wt{f}_{1}\\
&
=f_{1}^{2}+f_{2}^{2}+f_{1}(\theta)f_{2}-f_{2}(\theta)f_{1}-(c_{12}^{1}-f_{1}(\theta))f_{2}+(c_{12}^{2}-f_{2}(\theta))f_{1}\\
&=f_{1}^{2}+f_{2}^{2}-c_{12}^{1}f_{2}+c_{12}^{2}f_{1}
=\lapl_{f}.
\end{align*}
\end{proof}
Notice from \eqref{eq:lapl1} that the sum of squares is not an intrinsic operator of the sub-Riemannian structure.

\brem The same argument provides a proof of the fact that, on a 2-dimensional Riemannian manifold $M$ with local orthonormal frame $f_{1},f_{2}$ that satisfies
$$[f_{1},f_{2}]=a_{1}f_{1}+a_{2}f_{2},$$
the Laplace-Beltrami operator is locally expressed as
$$\lapl= f_{1}^{2}+f_{2}^{2}+a_{1}f_{2}-a_{2}f_{1}.$$
\erem

\section{The nilpotent approximation}\label{s:nilpapp}
In this section we briefly recall the concept of nilpotent approximation. For details one can see \cite{agrachevlocal,bellaiche,nostrolibro}.
\subsection{Privileged coordinates}
Let $\sub=(M,\distr,\metr)$ be a sub-Riemannian manifold and $(f_{1},\ldots,f_{k})$ an orthonormal frame. Fix a point $q\in M$ and consider the flag of the distribution 
$\distr^{1}_{q}\subset\distr^{2}_{q}\subset\ldots\subset\distr^{m}_{q}$. Recall that $k_i=\tx{dim}\, \distr^{i}_{q}-\tx{dim}\, \distr^{i-1}_{q}$ for $i=1,\ldots,m$, and that $k_{1}+\ldots+k_{m}=n$.

Let $O_{q}$ be an open neighborhood of the point $q\in M$. We say that a system of coordinates $\psi: O_{q}\to \R^{n}$ is \emph{linearly adapted} to the flag if, in these coordinates, we have $\psi(q)=0$ and
$$\psi_{*}(\distr^{i}_{q})=\R^{k_{1}}\oplus \ldots \oplus \R^{k_{i}}, \qquad \all i=1,\ldots,m.$$

Consider now the splitting
$\R^{n}=\R^{k_{1}}\oplus \ldots \oplus \R^{k_{m}}$
and denote its elements $x=(x_{1},\ldots,x_{m})$ where $x_{i}=(x_{i}^{1},\ldots,x_{i}^{k_{i}})\in \R^{k_{i}}$. 
The space of all differential operators in $\R^{n}$ with smooth coefficients forms an associative
algebra with composition of operators as multiplication. The differential operators with
polynomial coefficients form a subalgebra of this algebra with generators $1, x_{i}^{j} ,\frac{\partial}{\partial x_{i}^{j}},$ where
$i=1,\ldots,m;\  j = 1,\ldots,k_{i}$. We define weights of generators as
$$\nu(1)=0, \qquad \nu(x_{i}^{j})=i, \qquad \nu(\frac{\partial}{\partial x_{i}^{j}})=-i,$$ and the weight of monomials
$$\nu(y_{1}\cdots y_{\alpha}\frac{\partial^{\beta}}{\partial z_{1} \cdots \partial z_{\beta}})=\sum_{i=1}^{\alpha}\nu(y_{i})-\sum_{j=1}^{\beta}\nu(z_{j}).$$
Notice that a polynomial differential operator homogeneous with respect to $\nu$ (i.e. whose monomials are all of same weight) is homogeneous with respect to dilations $\delta_{t}:\R^{n}\to \R^{n}$ defined by
\bqn \label{dil} \delta_{t}(x_{1},\ldots,x_{m})=(tx_{1},t^{2}x_{2},\ldots, t^{m}x_{m}), \qquad t>0.
\eqn
In particular for a homogeneous vector field $X$ of weight $h$ it holds 
\bqn \label{eq:homogeneityvf}
\delta_{t*}X=t^{-h}X.
\eqn
A smooth vector field $X\in \tx{Vec}(\R^{n})$, as a first order differential operator, can be written as
$$X=\sum_{i,j} a_{i}^{j}(x) \frac{\partial}{\partial x_{i}^{j}},$$
and considering its Taylor expansion at the origin we can write the formal expansion
$$X\approx \sum_{h=-m}^{\infty} X^{(h)},$$
where $X^{(h)}$ is the homogeneous part of degree $h$ of $X$ (notice that every monomial of a first order differential operator has weight not smaller than $-m$).
Define the filtration of  $\tx{Vec}(\R^{n})$ $$\mc{D}^{(h)}=\{X\in \tx{Vec}(\R^{n}): X^{(i)}=0, \all i<h\}, \qquad \ell\in \mb{Z}.$$
\bdeff A system of coordinates $\psi: O_{q}\to \R^{n}$ defined near the point $q$ is said \emph{privileged} for a sub-Riemannian structure $\sub$ if these coordinates are linearly adapted to the flag and such that $\psi_{*}f_{i}\in \mc{D}^{(-1)}$ for every $i=1,\ldots,k$.
\edeff


Existence of privileged coordinates is proved, e.g. in \cite{agrachevlocal,
bellaiche}.
Notice however that privileged coordinates are not unique.

\bdeff Let $\sub=(M,\distr,\metr)$ be a regular sub-Riemannian manifold and $\{f_{1},\ldots,f_{k}\}$ a local orthonormal frame near a point $q$. Fixed a system of privileged coordinates, we define the \emph{nilpotent approximation of $\sub$ near $q$}, denoted by $\widehat{\sub}_{q}$, the sub-Riemannian structure on $\R^{n}$ having $\{\widehat{f}_{1}, \ldots, \widehat{f}_{k}\}$ as an orthonormal frame, where $\widehat{f}_{i}:=(\psi_{*}f_{i})^{(-1)}$.
\edeff
In what follows we often omit the coordinate map in the notation above, by writing $\widehat{f}_{i}:=f_{i}^{(-1)}$.
\brem \label{palla} It is well known that under the regularity hypothesis, $\widehat{\sub}_{q}$ is naturally endowed with a Lie group structure whose Lie algebra is generated by left-invariant vector fields $\widehat{f}_{1}, \ldots, \widehat{f}_{k}$. Moreover the sub-Riemannian distance $\widehat{d}$ in $\widehat{\sub}_{q}$ is homogeneous with respect to dilations $\delta_{t}$, i.e. $\widehat{d}(\delta_{t}(x),\delta_{t}(y))=t \,\widehat{d}(x,y)$. In particular, if $\widehat{B}_{q}(r)$ denotes the ball of radius $r$ in $\widehat{\sub}_{q}$, this implies $\delta_{t}(\widehat{B}_{q}(1))=\widehat{B}_{q}(t)$.
\erem

The following Lemma shows in which sense the nilpotent approximation is the first order approximation of the sub-Riemannian structure.
 
\bl\label{l:nilpotentvf} Let $M$ be a sub-Riemannian manifold and $X\in \VecM$. Fixed a system of privileged coordinates, we define $X^{\eps}:=\eps \delta_{\frac{1}{\eps}*}X$. Then
$$X^{\eps}= \wh{X}+\eps Y^{\eps}, \qquad \text{where} \ Y^{\eps} \ \text{is smooth w.r.t. } \eps.$$
\el 
\begin{proof} Since we work in a system of priviliged coordinates, in the homogeneous expansion of $X$ only terms of order $\geq 1$ appear. Hence we can write
$$X \simeq X^{(-1)}+X^{(0)}+X^{(1)}+\ldots$$ 
Applying the dilation  and using property \eqref{eq:homogeneityvf} we get
$$\delta_{\frac{1}{\eps}*}X \simeq \frac{1}{\eps}X^{(-1)}+X^{(0)}+\eps X^{(1)}+\ldots$$ 
Multiplying by $\eps$ and using that, by definition, $\wh{X}=X^{(-1)}$, we have
\bqn \label{eq:feps}
X^{\eps} \simeq \wh{X}+\eps X^{(0)}+\eps^{2} X^{(1)}+\ldots
\eqn

\end{proof}

In other words the nilpotent approximation of a vector field at a point $q$ is the first meaningful term that appears in the expansion when one consider the blow up coordinates near the point $q$, with rescaled distances. 

\brem \label{rem:heis}
If $\sub=(M,\distr,\metr)$ is a 3D contact sub-Riemannian manifold, then $\dim \distr_{q}=2$ and $\dim \distr^{2}_{q}=3$ for all $q\in M$. Under this assumption the nilpotent approximation $\wh{\sub}_{q}$ of $\sub$ at every point $q\in M$ is isometric to the Heisenberg group, since this is the only nilpotent left-invariant structure with $\mc{G}(\sub)=(2,3)$ (see e.g. \cite{miosr3d} for a classification of left-invariant structures on 3D Lie groups). 

The sub-Riemannian structure on the Heisenberg group is defined by the global orthonormal frame on $\R^{3}$ (where coordinates are denoted by $q=(x,y,w)$) 
\begin{align} \label{eq:heisenberg}
\wh{f}_{1}&=\partial_{x}+\frac{y}{2}\partial_{w}, \qquad \wh{f}_{2}=\partial_{y}-\frac{x}{2}\partial_{w}.
\end{align} 
\erem
Notice that the Lie algebra $\tx{Lie}\{\wh{f}_{1},\wh{f}_{2}\}$ is nilpotent since $[\wh{f}_{1},[\wh{f}_{1},\wh{f}_{2}]]=[\wh{f}_{2},[\wh{f}_{1},\wh{f}_{2}]]=0$. Moreover the Reeb vector field is $\wh{f}_{0}=\partial_{w}$ and the local invariants of the structure are identically zero $\chi=\kappa=0$.

\subsection{Normal coordinates} 
In the 3D contact case there exists a smooth normal form of the sub-Riemannian structure (i.e. of its orthonormal frame) which is the analogue of normal coordinates in Riemannian geometry. This normal form is crucial for the study of  the heat kernel of the sub-Laplacian with a perturbative approach, since it presents the sub-Riemannian structure of a general 3D contact case as a perturbation of the Heisenberg (nilpotent) case. 

\bt[\hspace{-0.12cm} \cite{srR3, dido}] \label{t:normal}
Let $M$ be a 3D contact sub-Riemannian manifold and $f_{1},f_{2}$ a local orthonormal frame. There exists a smooth coordinate system $(x,y,w)$ such that
\begin{align*}
f_{1}&=(\partial_{x}+\frac{y}{2}\partial_{w})+\beta y (y\partial_{x}-x\partial_{y})+\g y \partial_{w},\\[0.1cm]
f_{2}&=(\partial_{y}-\frac{x}{2}\partial_{w})-\beta x (y\partial_{x}-x\partial_{y})+\g x \partial_{w},
\end{align*}
where $\beta=\beta(x,y,w)$ and $\gamma=\gamma(x,y,w)$ are smooth functions that satisfy the following boundary conditions
$$\beta(0,0,w)=\g(0,0,w)=\frac{\partial \g}{\partial x}(0,0,w)=\frac{\partial \g}{\partial y}(0,0,w)=0.$$ 
\et

Notice that the normal coordinate system 
is privileged at 0. 
Indeed from the explicit expression of the frame
it immediately follows that these coordinates are linearly adapted at $0$ since 
$$\distr_{0}=\tx{span}\{f_{1}(0),f_{2}(0)\}=\tx{span}\{\partial_{x},\partial_{y}\}=\R^{2}.$$
Moreover the weights of the coordinates $(x,y,w)$ at the origin are
$$\nu(x)=\nu(y)=1, \qquad \nu(w)=2,$$
and every homogeneous term of the vector fields $f_{1},f_{2}$ has degree $\geq -1$.

Finally, notice also that, when $\beta=\gamma=0$,  we recover the Heisenberg group structure \eqref{eq:heisenberg}.

\section{Perturbative method} \label{sez:pert}
In this section we consider the sub-Riemannian heat equation associated to a sub-Riemannian structure $f$ on a complete sub-Riemannian manifold $M$, i.e. the initial value problem
\bqn \label{eq:subheat}
\begin{cases}
\dfrac{\partial \psi}{\partial t}(t,x)= \lapl_{f} \psi(t,x), \qquad \text{in} \   (0,\infty) \times M,\\[0.3cm]
\psi(0,x)=\ph(x), \qquad \qquad \   x\in M, \quad \ph\in C^{\infty}_{0}(M).
\end{cases}
\eqn
where $\psi(0,x)=\underset{t \to 0}{\lim}\, \psi(t,x)$ and the limit is meant in the distributional sense.

\medskip
Recall that a differential operator $\mc{L}$ is said to be hypoelliptic on a subset $U\subset M$ if every distributional solution to $\mc{L}u=\phi$ is $C^{\infty}(U)$, whenever $\phi \in C^{\infty}(U)$. 
The following well-known H\"ormander Theorem  gives a sufficient condition for the hypoellipticity of a second order differential operator.


\bt[H\"ormander,\cite{hormander}] Let $\mc{L}$ be a differential operator on a manifold $M$, that locally in a neighborhood $U$ is written as $$\mc{L}=\sum_{i=1}^{k} X_{i}^{2}+X_{0},$$ where $X_{0},X_{1},\ldots,X_{k}\in \VecM$. If $\tx{Lie}_{q}\{X_{0},X_{1},\ldots,X_{k}\} = T_{q}M$ for all $q \in U$, then $\mc{L}$ is hypoelliptic.
\et\
From this Theorem and the bracket generating condition it follows that the sub-laplacian $\lapl_{f}$ is hypoelliptic. Moreover, since $M$ is complete and the sub-Laplacian is symmetric and negative with respect to the Popp's measure (see \eqref{eq:sym}), it follows that $\lapl_{f}$ is essentially self-adjoint on $C^{\infty}_{0}(M)$ (see also (iii) in Remark \ref{r:superem}).

As a consequence the operator $\lapl_{f}$ admits a unique self-adjoint extension on $L^{2}(M)$ and the heat semigroup $\{e^{t\lapl_{f}}\}_{t\geq 0}$ is a well-defined one parametric family of bounded operators on $L^{2}(M)$. Moreover the heat semigroup is contractive on $L^{2}(M)$ (see \cite{strichartz}).

The problem \eqref{eq:subheat} has a unique solution,  for every initial datum $\ph\in L^{2}(M)$, namely $\psi(t,x):=e^{t\lapl_{f}}\ph$. Due to the hypoellipticity of $\lapl_{f}$, the function $(t,x) \mapsto e^{t\lapl_{f}}\ph(x)$ is smooth on $(0,\infty) \times M$ and

$$e^{t\lapl_{f}}\ph(x) = \int_{M} p(t,x,y)\ph(y)dy,	\qquad \ph \in C^{\infty}_{0}(M),$$
where $p(t,x, y)$ is the so-called \emph{heat kernel} associated to $e^{t\lapl_{f}}$, that satisfies the following properties 
\bi
\iii[(i)] $p(t,x,y) \in C^{\infty}(\R^{+}\times M\times M)$,
\iii[(ii)] $p(t,x,y)=p(t,y,x), \quad \all t>0, \all x,y\in M$,
\iii[(iii)] $p(t,x,y)>0, \quad \all t>0, \all x,y\in M$.
\ei
A probabilistic approach to hypoelliptic diffusion equation can be found in \cite{baudoin,bismut,kusustroock}, where the existence of a smooth heat kernel for such equations is given.

In particular, since $e^{t\lapl_{f}}\ph$ satisfies the initial condition, it holds
$$\lim_{t\to 0} \int_{M} p(t,x,y) \ph(y) dy=\ph(x).$$
and the heat kernel $p(t,x,y)$ is a solution of the problem \eqref{eq:subheat} with the initial condition $\psi(0,x)=\delta_{y}(x)$, where $\delta_{y}$ denotes the Dirac delta function.
%

For a more detailed discussion on the analytical properties of the sub-Riemannian heat equation and its heat kernel one can see  \cite{strichartz, garofalob}.  

To study the asymptotics of the heat kernel associated to the sub-Riemannian structure defined by $f$ near a point $q\in M$, we consider the approximation of the sub-Riemannian structure (cfr. also Lemma \ref{l:nilpotentvf}).
\bdeff 
Let $f_{1},\ldots,f_{k}$ be an orthonormal frame for a sub-Riemannian structure on $M$ and fix a system of privileged coordinates around the point $q\in M$. The $\eps$-approximated system at $q$ is the sub-Riemannian structure induced by the orthonormal frame 
$f_{1}^{\eps},\ldots,f_{k}^{\eps}$ defined by
$$f_{i}^{\eps}:=\eps \delta_{1/\eps*}f_{i}, \qquad 
i=1,\ldots,k.$$
\edeff

\brem\label{r:f0}
Notice that  in the definition of approximated structure we have to perturbate the basis of the distribution only. Their Lie brackets are changed accordingly to the formula
$$[f_{i}^{\eps},f_{j}^{\eps}]=\eps^{2}\delta_{1/\eps*}[f_{i},f_{j}].$$
In particular, in the 3D contact case, the Reeb vector field $f_{0}^{\eps}$ of the $\eps$-approximated structure is related to the unperturbed one by $f_{0}^{\eps}=\eps^{2}\delta_{1/\eps*}f_{0}$.
\erem
The following Lemma
shows the relation between the heat kernel of $f$ and the one defined by $f^{\eps}$, whose existence is guaranteed by the above results. 

\bl Let $M$ be a sub-Riemannian manifold and fix a set of privileged coordinates in a neighborhood $N$ of $q$. Denote by $f$ the sub-Riemannian structure and by $f^{\eps}$ its $\eps$-approximation at $q$.
If we denote by $p(t,x,y)$ and $p^{\eps}(t,x,y)$ the heat kernels respectively of the equations
$$\frac{\partial \psi}{\partial t}(t,x)= \lapl_{f} \psi(t,x), \qquad \frac{\partial \psi}{\partial t}(t,x)= \lapl_{f^{\eps}} \psi(t,x),$$
we have that
$$p^{\eps}(t,x,y)=\eps^{Q}p(\eps^{2}t,\delta_{\eps}x,\delta_{\eps}y), \qquad \all x,y \in N,$$
where $Q$ denotes the Hausdorff dimension of $M$.
\el
\begin{proof}
Recall that, if $f_{1},\ldots,f_{k}$ is an orthonormal frame for the sub-Riemannian structure, the approximated system is defined by vector fields $$f_{i}^{\eps}:=\eps \delta_{1/\eps*}f_{i}, \qquad 
i=1,\ldots,k.$$
Then the proof follows from the following facts

$(i)$. If we perform the change of coordinates $x^{\prime}=
\dil_{1/\eps}x$ and we denote by $q^{\eps}(t,x,y)$ the heat kernel written in the new coordinate system (which depends on $\eps$), we have the equality
$$q^{\eps}(t,x,y)=\eps^{Q}p(t,\dil_{\eps}x,\dil_{\eps}y)$$ 
Indeed it is easy to see that $q^{\eps}(t,x,y)$ is a solution of the equation written in the new variables $x^{\prime}=\delta_{1/\eps}x$ and the factor $\eps^{Q}$ comes from the fact that $|\det \dil_{\eps}|=\eps^{Q}$. More precisely if we set
$$\psi^{\eps}(t,x)=\int_{M} q^{\eps}(t,x,y) \ph^{\eps}(y)dy,$$ 
where $\ph^{\eps}(x)=\ph(\delta_{\eps}x)$, one can see that the initial condition is satisfied:
\begin{align*}
\lim_{t\to 0} \int_{M} q^{\eps}(t,x,y) \ph^{\eps}(y) dy &= \lim_{t\to 0} \int_{M} \eps^{Q}p(t,\dil_{\eps}x,\dil_{\eps}y) \ph(\dil_{\eps}y) dy\\
&=\lim_{t\to 0} \int_{M} p(t,\dil_{\eps}x,z) \ph(z) dz \qquad(z=\delta_{\eps}y)\qquad\\
&= \ph(\dil_{\eps}x)=\ph^{\eps}(x).
\end{align*}

$(ii)$. 
Since $\delta_{\frac{1}{\eps}*}f_{i}= \frac{1}{\eps}f_{i}^{\eps}$, for every $i=1,\ldots,k$, using \eqref{eq:hom2} we get
$$\frac{\partial \psi}{\partial t}= \lapl_{\frac{1}{\eps}f^{\eps}} \psi=\frac{1}{\eps^{2}}\Delta_{f^{\eps}} \psi.$$
This equality can be rewritten as
$$\eps^{2}\frac{\partial \psi}{\partial t}= \lapl_{f^{\eps}} \psi,$$
and performing the change of variable $t=\eps^{2}\tau$ we get
$$\frac{\partial \psi}{\partial \tau}= \lapl_{f^{\eps}} \psi,$$
from which it follows
$$ p^{\eps}(t,x,y)=q^{\eps}(\eps^{2}t,x,y)=\eps^{Q}p(\eps^{2}t,\delta_{\eps}x,\delta_{\eps}y).$$


%




\end{proof}

This result is useful for the study of the asymptotics on the diagonal since, by construction, the intial point is fixed for the dilation $\delta_{\eps}$. In particular we can recover the small time behaviour of the original heat kernel from the approximated one
\bc \label{c:cc} Following the notations introduced above
$$p^{\eps}(1,0,0)=\eps^{Q}p(\eps^{2},0,0).$$
\ec

As a corollary we recover an homogeneity property of the heat kernel for a nilpotent structure. 
 \bc Assume that the sub-Riemannian structure is regular and nilpotent. Then the heat kernel $p(t,x,y)$ of the heat equation satisfies the homogeneity property
$$p(t,x,y)=\lam^{Q}p(\lam^{2}t, \delta_{\lam}x,\delta_{\lam}y), \qquad \all \lam >0.$$
\ec
This result was already well-known (see e.g. \cite{folland}).
In our notations, it is a consequence of the fact that, if $f=\wh{f}$, then $f=f^{\eps}=\wh{f}$ for every $\eps>0$.

\subsection{General method} \label{sez:pert0}
In this section we briefly recall the perturbative method for the heat equation presented in Chapter 3 of \cite{rosenberg}. For further discussions one can see also \cite{spectrebook,cycon} and \cite{yu}. 

Let $X,Y$ be operators on a Hilbert space of functions. Moreover  assume that $X$ and $X+Y$ (where $Y$ is treated as a perturbation of $X$) have well defined heat operators $e^{tX}$, $e^{t(X+Y)}$, i.e. a semigroup of one parameter family of bounded self-adjoint operators satisfying
$$(\partial_{t}-X)e^{tX}\ph=0, \qquad \lim_{t\to 0}e^{tX}\ph=\ph,$$
and similarly for $X+Y$. 
Given $A(t),B(t)$ two operators on the Hilbert space, if we denote their convolution as
$$(A*B)(t)=\int_{0}^{t}A(t-s)B(s)ds,$$
then the classical Duhamel formula 
$$e^{t(X+Y)}=e^{tX}+\int_{0}^{t}e^{(t-s)(X+Y)}Ye^{sX} ds,$$
can be rewritten as follows
$$e^{t(X+Y)}=e^{tX}+e^{t(X+Y)}*Ye^{tX}.$$
Iterating this construction one gets the expansion
\bqn \label{duhamel}e^{t(X+Y)}=e^{tX}+e^{tX}*Ye^{tX}+e^{t(X+Y)}*(Ye^{tX})^{*2}.\eqn
where $A^{*2}=A*A$ denotes the iterated convolution product.

\vspace{0.1cm}
If $A(t)$ and $B(t)$ have heat kernels $a(t,x,y)$ and $b(t,x,y)$ respectively, then $(A*B)(t)$ has kernel (see again \cite{rosenberg,yu})
$$(a*b)(t,x,y):= \int_{0}^{t} \int_{M} a(s,x,z)b(t-s,z,y) dzds,$$
Interpreting \eqref{duhamel} at the level of kernels, denoting by $p(t,x,y)$ the heat kernel for the operator $X$ and by $p^{Y}(t,x,y)$ the kernel of the perturbed operator $X+Y$ we can write the expansion
\begin{align} \label{p:duhamel} p^{Y}(t,x,y)
&=p(t,x,y)+(p*Yp)(t,x,y)+ (p^{Y}*Yp*Yp)(t,x,y)
\end{align}


%



\section{Proof of Theorem \ref{t:main}}
In this section we compute the first terms of the small time asymptotics of the heat kernel. To this extent we compute the sub-Laplacian associated to the approximated sub-Riemannian structure and we use the pertubative method of Section \ref{sez:pert0} to compute this terms using the explicit expression of the heat kernel in the Heisenberg group.

\brem
The sub-Laplacian on the Heisenberg group $H^{3}$ is written as the sum of squares (cfr. also Remarks \ref{r:superem} and \ref{rem:heis})
$$\lapl_{\wh{f}}=\wh{f}_{1}^{2}+\wh{f}_{2}^{2}=(\partial_{x}-\frac{y}{2}\partial_{w})^{2}+ (\partial_{y}+\frac{x}{2}\partial_{w})^{2}. $$
The heat kernel for $\Delta_{\wh{f}}$ has been computed explicitly for the first time in \cite{gaveau}. Here we use the expression given in \cite{laplacian} in the same coordinate set. Denote by $q=(x,y,w)\in \R^{3}$ a point in the Heisenberg group.
The heat kernel $H(t,q,q^{\prime})$, is presented as 
\bqn \label{eq:Hh} H(t,q,q^{\prime})=h_{t}(q^{\prime} \circ q^{-1}),  \eqn
where
\bqn \label{eq:explheis}
h_{t}(x,y,w)=\frac{1}{2(2\pi t)^2}\int_\R \frac{s}{\sinh s}\exp\left({-\frac{s(x^2+y^2)}{4t\tanh s}}\right)\cos(\frac{ws}{t}) ds,
\eqn
and $\circ$ denotes the group law in $H_{3}$
\bqn \nn
(x,y,w)\circ(x^{\prime},y^{\prime},w^{\prime})=(x+x^{\prime},y+y^{\prime},w+w^{\prime}+\frac{1}{2}(x^{\prime}y-xy^{\prime})).
\eqn
Notice that the inverse of an element with respect to $\circ$ is 
\bqn \nn
(x,y,w)^{-1}=(-x,-y,-w).
\eqn

For a discussion on the convergence of the integral \eqref{eq:explheis} one can see \cite{bealsgaveau}.
\erem

\subsection{Local invariants}
In this section we compute the invariants $\chi$ and $\kappa$ at the origin of the sub-Riemannian manifold. By Theorem \ref{t:normal} we can assume that the orthonormal frame has the form
\begin{align} 
f_{1}&=(\partial_{x}+\frac{y}{2}\partial_{w})+\beta y (y\partial_{x}-x\partial_{y})+\g y \partial_{w},\nn \\
f_{2}&=(\partial_{y}-\frac{x}{2}\partial_{w})-\beta x (y\partial_{x}-x\partial_{y})+\g x \partial_{w}, \label{eq:ultima}
\end{align}
where $\beta$ and $\gamma$ are smooth functions near $(0,0,0)$ that satisfy
\bqn \label{eq:bound}
\beta(0,0,w)=\g(0,0,w)=\frac{\partial \g}{\partial x}(0,0,w)=\frac{\partial \g}{\partial y}(0,0,w)=0.
\eqn

Moreover, since we are interested up to second order terms in the expansion \eqref{eq:feps}, we can assume the following
\bl \label{lemmau}We can assume that the orthonormal frame has the form
\begin{align}
f_{1}&=\partial_{x}-\frac{y}{2}(1+\g)\partial_{w}, \nn \\
f_{2}&=\partial_{y}+\frac{x}{2}(1+\g)\partial_{w}, \label{eq:normalred}
\end{align}
where $\g$ is a quadratic polynomial of the form $\g(x,y)=ax^{2}+bxy+cy^{2}$, for some $a,b,c\in \R$.
\el
\begin{proof} Recall that if we expand a vector field $X$ in homogeneous components (when written in a privileged coordinate system)
$$X \simeq X^{(-1)}+X^{(0)}+X^{(1)}+X^{(2)}+\ldots$$ 
its $\eps$-approximation $X^{\eps}$ has the following expansion
\bqn \label{eq:feps}
X^{\eps} \simeq \wh{X}+\eps X^{(0)}+\eps^{2} X^{(1)}+\eps^{3}X^{(2)}+\ldots
\eqn

It is then sufficient to consider, in the Taylor expansion of the orthonormal frame near the origin, only the homogeneous term up to weight one, since every other term with weight $\geq 2$ gives a contribution $o(\eps^{2})$ when one compute the expansion of $\lapl_{f^{\eps}}$. Hence a contribution $o(\eps^{2})$ in the heat kernel due to \eqref{p:duhamel}.
 
Moreover the boundary condition \eqref{eq:bound} implies that the following derivatives of the coefficients of \eqref{eq:ultima} vanish at the origin
\bqn \label{eq:bound1}
\frac{\partial \beta}{\partial w}(0,0,0)=\frac{\partial \g}{\partial w}(0,0,0)=\frac{\partial^{2} \g}{\partial w \partial x}(0,0,0)=\frac{\partial^{2} \g}{\partial w \partial y}(0,0,0)=0,
\eqn
together with all higher order derivatives with respect to $w$.

Since $\nu(\partial_{x})=\nu(\partial_{y})=1$ and $\beta(0,0,0)=0$, the terms $\beta y (y\partial_{x}-x\partial_{y})$ and $\beta x (y\partial_{x}-x\partial_{y})$ have weight $\geq 2$. Moreover $\nu(\partial_{w})=-2$ implies that $x\partial_{w}$ and $ y \partial_{w}$ have weight $-1$. The only terms that we need in the expansion of $\g$ are those of weight less or equal than one. Since the terms of order zero vanish by \eqref{eq:bound} and \eqref{eq:bound1}, the only meaningful term in the expansion of $\g$ is
\begin{align*}
\g(x,y,w) \sim&\ 
\frac{\partial^{2} \g}{\partial x^{2}}x^{2}+\frac{\partial^{2} \g}{\partial x y}xy+\frac{\partial^{2} \g}{\partial y^{2}}y^{2}
\end{align*}
where derivatives are computed at the origin (0,0,0).
\end{proof}

Now we express the invariants $\chi$ and $\kappa$ in terms of the perturbation \eqref{eq:normalred}.
\bl Assume that the orthonormal frame of the sub-Riemannian structure has the form \eqref{eq:normalred}. Then value of the invariants at the origin are
\bqn
\chi=2\sqrt{b^{2}+(c-a)^{2}}, \qquad 
\kappa=2(a+c).
\eqn
\el
\begin{proof}
To compute the invariants we need to compute the Reeb vector field $f_{0}$ and the structure constant of the Lie algebra $\tx{Lie}\{f_{0},f_{1},f_{2}\}$. Every contact form for the structure is a multiple of 
$$\wt{\omega}=dz-\frac{x}{2}(1+\g)dy +\frac{y}{2}(1+\g)dx , $$
whose differential is computed as follows
\begin{align*}
d\wt{\omega}=-(1+2\g) dxdy.
\end{align*}
Since 
$d\wt{\omega}(f_{1},f_{2})=-(1+2\g)$,
the normalized contact form $\omega$ that satisfies $d\omega(f_{1},f_{2})=1$ (see Remark \ref{rem:contact} and \eqref{eq:setting}) is
\begin{align*}\omega:=-\frac{1}{1+2\g}\wt{\omega}=-\frac{1}{1+2\g}(dz-\frac{x}{2}(1+\g)dy +\frac{y}{2}(1+\g)dx).
\end{align*}
Notice that every contact form vanishes on the distribution. Thus $$d(\phi \omega)(f_{1},f_{2})=\phi\,d\omega(f_{1},f_{2}), \qquad \all \phi\in C^{\infty}(M).$$
Next we compute the differential of the normalized contact form
\begin{align*}
d\omega&=\frac{d(1+2\g)\wedge \til{\omega}}{(1+2\g)^{2}}+dxdy\\
&=\frac{2 \partial_{x}\g}{(1+2\g)^{2}}dxdw+\frac{2\partial_{y}\g}{(1+2\g)^{2}}dydw+\left(1-\frac{2\g(1+\g)}{(1+2\g)^{2}}\right)dxdy,
\end{align*} 
The Reeb vector field is, by definition, the kernel of $d\omega$ (normailzed in such a way that $\omega(f_{0})=1$). From this one gets
$$f_{0}=\frac{2 \partial_{x}\g}{1+2\g}\partial_{y}-\frac{2\partial_{y}\g}{1+2\g}\partial_{x}+\left(\frac{2\g(1+\g)}{1+2\g}-(1+2\g)\right)\partial_{w}.$$
The commutator between horizontal vector fields is computed as follows
\begin{align*}
[f_{2},f_{1}]&=\left[\partial_{y}+\frac{x}{2}(1+\g)\partial_{w},\partial_{x}-\frac{y}{2}(1+\g)\partial_{w}\right]\\
&=-(1+2\g)\partial_{w},
\end{align*}
and writing $[f_{2},f_{1}]=f_{0}+c_{12}^{1}f_{1}+c_{12}^{2}f_{2}$ we find the structure constants
$$c_{12}^{1}=\frac{2\partial_{y}\g}{1+2\g},\qquad c_{12}^{2}=-\frac{2\partial_{x}\g}{1+2\g}.$$

Moreover, a longer computation for $[f_{1},f_{0}]$ and $[f_{2},f_{0}]$ shows that
 $$c_{0i}^{j}=-\frac{2}{(1+2\g)^{2}}\til{c}_{0i}^{j},$$ where we set
\begin{align*}
\til{c}_{01}^{1}&=(1+2\g)\partial_{xy}\g-2\partial_{y}\g\partial_{x}\g,\\
\til{c}_{01}^{2}&=-(1+2\g)\partial_{xx}\g+2(\partial_{x}\g)^{2},\\
\til{c}_{02}^{1}&=(1+2\g)\partial_{yy}\g-2(\partial_{y}\g)^{2},\\
\til{c}_{02}^{2}&=-(1+2\g)\partial_{xy}\g+2\partial_{y}\g\partial_{x}\g.
\end{align*}
Recalling that at the origin $\partial_{x}\g=\partial_{y}\g=0$, while $\partial_{xx}\g=2a, \partial_{xy}\g=b, \partial_{yy}\g=2c$
it follows from \eqref{eq:chi} that
$$\chi=2\sqrt{-\det
\begin{pmatrix}
b& c-a\\
c-a &-b
\end{pmatrix}}=2\sqrt{b^{2}+(c-a)^{2}}.
$$
and
\begin{align*}
\kappa&=f_{2}(c_{12}^{1})-f_{1}(c_{12}^{2})-(c_{12}^{1})^{2}-(c_{12}^{2})^{2}+\frac{c_{01}^{2}-c_{02}^{1}}{2}=2(a+c).
\end{align*}
\end{proof}
\subsection{Asymptotics}
In this section we compute the Laplacian $\lapl_{f^{\eps}}$ up to second order in $\eps$.
First notice that
\begin{align*}
f_{1}^{\eps}&=(\partial_{x}-\frac{y}{2}\partial_{w})-\eps^{2}(\frac{y}{2}\g\partial_{w})+o(\eps^{2}),\\
f_{2}^{\eps}&=(\partial_{y}+\frac{x}{2}\partial_{w})+\eps^{2}(\frac{x}{2}\g\partial_{w})+o(\eps^{2}).
\end{align*}
Moreover, defining
$f_{0}^{\eps}:=\eps^{2} \delta_{\frac{1}{\eps}*} f_{0}$ (see Remark \ref{r:f0}), from the formula 
$$[f_{2}^{\eps},f_{1}^{\eps}]=f_{0}^{\eps}+(c_{12}^{1})^{\eps}f_{1}^{\eps}+(c_{12}^{2})^{\eps}f_{2}^{\eps},$$
we get the following expansion 
\begin{align}
(c_{12}^{1})^{\eps}&=\eps \frac{2\eps (\partial_{y}\g)}{1+2\eps^{2}\g}=2\eps^{2}\partial_{y}\g+o(\eps^{2}), \label{eq:cije1}\\
(c_{12}^{2})^{\eps}&=-\eps \frac{2\eps (\partial_{x}\g)}{1+2\eps^{2}\g}=-2\eps^{2}\partial_{x}\g+o(\eps^{2}). \label{eq:cije2}
\end{align}

Thus we can compute every term defining the sub-Laplacian
\begin{align*}
(f_{1}^{\eps})^{2}&=(\partial_{x}-\frac{y}{2}\partial_{w})^{2}-\eps^{2} \left( (\partial_{x}-\frac{y}{2}\partial_{w})(\frac{y}{2}\g\partial_{w})+(\frac{y}{2}\g\partial_{w})(\partial_{x}-\frac{y}{2}\partial_{w})\right)+o(\eps^{4})\\
&=(\wh{f}_{1})^{2}-\eps^{2}(y\g \partial_{wx}-\frac{y^{2}}{2}\g \partial_{w}^{2}+\frac{y}{2}\partial_{x}\g \partial_{w})+o(\eps^{4}),\\[0.4cm]
(f_{2}^{\eps})^{2}&=(\partial_{y}+\frac{x}{2}\partial_{w})^{2}+\eps^{2} \left( (\partial_{y}+\frac{x}{2}\partial_{w})(\frac{x}{2}\g\partial_{w})+(\frac{x}{2}\g\partial_{w})(\partial_{y}+\frac{x}{2}\partial_{w})\right)+o(\eps^{4})\\
&=(\wh{f}_{2})^{2}+\eps^{2}(x\g \partial_{wx}+\frac{x^{2}}{2}\g \partial_{w}^{2}+\frac{x}{2}\partial_{y}\g \partial_{w})+o(\eps^{4}),
\end{align*}
where $\wh{f}_{1}, \wh{f}_{2}$ denote the orthonormal frame of the Heisenberg group (see again Remark \ref{rem:heis}). Moreover, from \eqref{eq:cije1} and \eqref{eq:cije2} one easily gets
\begin{align*}
-(c_{12}^{1})^{\eps} f_{2}^{\eps}&=-2\eps^{2}\partial_{y}\g(\partial_{y}+\frac{x}{2}\partial_{w})+o(\eps^{4}),\\[0.3cm]
(c_{12}^{2})^{\eps} f_{1}^{\eps}&=-2\eps^{2}\partial_{x}\g(\partial_{x}-\frac{y}{2}\partial_{w})+o(\eps^{4}).
\end{align*}
Recollecting all the terms we find that
$$\lapl_{f^{\eps}}=\lapl_{\wh{f}}+\eps^{2}\Y+o(\eps^{2}),$$
where
$$\lapl_{\wh{f}}=\wh{f}_{1}^{2}+\wh{f}_{2}^{2}=(\partial_{x}-\frac{y}{2}\partial_{w})^{2}+ (\partial_{y}+\frac{x}{2}\partial_{w})^{2}, $$
is the sub-Laplacian on the Heisenberg group, and
$\Y$ denotes the second order differential operator
$$\Y=\frac{\g}{2}(x^{2}+y^{2})\partial_{w}^{2}+\g(x \partial_{wy}-y\partial_{wx})-\frac{1}{2}(x\partial_{y}\g-y\partial_{x}\g)\partial_{w}-2(\partial_{x}\g \partial_{x}+\partial_{y}\g \partial_{y}).$$

Specifying expansion \eqref{p:duhamel} to our case, where
$Y=Y(\eps)$ is a smooth perturbation which expands with respect to $\eps$ as follows 
$$Y(\eps)=\eps^{2}\Y+o(\eps^{2})$$
we find that
\bqn\label{45}
p^{\eps}=H+\eps^{2} (H*\Y H)+\eps^{4}(p^{\eps}*\Y H * \Y H).
\eqn

Since $p^{\eps}$ is the heat kernel of a contraction semigroup (for every $\eps>0$), one can see that the last term in \eqref{45} is bounded and
 $$p^{\eps}(1,0,0)=H(1,0,0)+\eps^{2}H*\Y H(1,0,0)+O(\eps^{4}),$$
 where $$\Y H(t,q,q^{\prime})=\Y_{q}H(t,q,q^{\prime}),$$ means that $\Y$ acts as a differential operator on the first spatial variable.

From the explicit expression \eqref{eq:explheis} it immediately follows that
$$H(1,0,0)= \frac{1}{16t^{2}}.$$

Thus, denoting by $K_{1}:=H*\Y H(1,0,0)$ from Corollary \ref{c:cc} we have the expansion of the original heat kernel
$$p(t,x,x)\sim \frac{1}{16 t^{2}}(1+ K_{1}t +O(t^{2})).$$

We are left to computation of the convolution between $H$ and $YH$, namely
\begin{align}
H*\Y H(1,0,0)&= \int_{0}^{t} \int_{\R^{3}} H(s,0,q)\Y H(1-s,q,0) \,dqds \nn \\
&= \int_{0}^{t} \int_{\R^{3}} h_{s}(q)\Y h_{1-s}(q) \,dqds, \qquad q=(x,y,w). \label{eq:ausilio}
\end{align}

Computing derivatives under the integral sign one gets
\begin{align*}
\Y h_{t}(x,y,w)=&-\frac{1}{(4 \pi  t)^2}\int_{\R}  \frac{r}{
 \sinh r} \exp\left({-\frac{r(x^2+y^2)}{2t\tanh r}}\right)  \frac{r}{t^2}\\
 &\times  \left[  \g(x,y) \cos(\frac{ r w}{t})\left( r (x^2 + y^2) -\frac{4t}{\tanh r}\right) + 
   t\, \g^{\prime}(x,y) \sin(\frac{r w}{t})\right] dr,
   \end{align*}
   
   where
   $$\g(x,y)=ax^{2}+bxy+cy^{2}, \qquad \g^{\prime}(x,y)=2(a-c)xy+b(y^{2}-x^{2}).$$   
   
Notice that we are interested only in computing how the integral \eqref{eq:ausilio} depends on the constants $a,b,c$, and the perturbation of the metric is on the variables $x,y$ only.  Hence, using that the integrand has exponential decay with respect to $x,y$, we can exchange the order of integration in \eqref{eq:ausilio} and integrate first with respect to these variables.

Using that, for $\al>0$
\begin{gather*}
\iint_{\R^{2}} x^{2} e^{-\al(x^{2}+y^{2})}dxdy=\frac{c_{1}\pi}{\al}, \qquad \iint_{\R^{2}} x^{4} e^{-\al(x^{2}+y^{2})}=\frac{c_{2}\pi}{\al^{3}},\\
\iint_{\R^{2}} xy e^{-\al(x^{2}+y^{2})}dxdy=0,
\end{gather*}
for some constants $c_{1},c_{2}>0$, it is easily seen that integrating \eqref{eq:ausilio} we get an expression of the kind
\bqn \label{eq:ausz}
H*\Y H(1,0,0)=C_{0} (a+c)= C_{0}\kappa,
\eqn
where $C_{0}$ is a universal constant that does not depend on the sub-Riemannian structure. Hence the value of $C_{0}$ can be computed from some explicit formula of the heat kernel in the non nilpotent case.
Using the expression given in \cite{bonnefontsu2} for the heat kernel on $SU(2)$, where the value of local invariants are constant $\chi=0, \kappa=1$ (see also \cite{miosr3d}) we get that
$$p_{SU(2)}(t,0,0)=\frac{e^{t}}{16t^{2}}\sim \frac{1}{16t^{2}}(1+t+O(t^{2})),$$
where we renormalized the constants in order to fit into our setting. Hence $C_{0}=1$ in \eqref{eq:ausz}, and the Theorem is proved.

\brem
The same method applies to get a quick proof of the following-well known result (see e.g. \cite{spectrebook,rosenberg}): on a 2-dimensional Riemannian manifold $M$, the heat kernel $p(t,x,y)$ satisfies an asymptotic expansion on the diagonal
$$p(t,x,x)\sim \frac{1}{4\pi t}(1+ \frac{K(x)}{6}t +O(t^{2})),\qquad \text{for } t \to 0,$$
where $K(x)$ denotes the gaussian curvature at the point $x\in M$. 
Indeed one can use the normal coordinates on $M$ to write the orthonormal frame in the following way
\begin{align*}
f_{1}&=\partial_{x}+\beta y (y\partial_{x}-x\partial_{y}), \\
f_{2}&=\partial_{y}-\beta x (y\partial_{x}-x\partial_{y}),
\end{align*}
where $\beta$ is, a priori, a smooth function $\beta=\beta(x,y)$. Reasoning as in Lemma \ref{lemmau}, $\beta$ can be chosen as a constant since we are interested only in first order term. In this case it is also easily seen that the Gaussian curvature at the origin is computed via the parameter $\beta$ as $K=6\beta$.

\erem


%

{\bf Aknowledgements}. The author is grateful to Andrei Agrachev for suggesting the problem and many helpful discussions.

{\small
\bibliography{bibliogr}
\bibliographystyle{abbrv}
}

\end{document}